\documentclass{llncs}

\usepackage[utf8]{inputenc}

\usepackage[disable,textsize=footnotesize]{todonotes}   

\usepackage{booktabs}
\usepackage{graphicx}
\usepackage[caption=false]{subfig}
\usepackage{tabularx, booktabs, multirow}
\usepackage{siunitx}
\usepackage{verbatim}
\newcommand{\figref}[1]{Fig.~\ref{#1}}
\newcommand{\eqreffull}[1]{\ensuremath{\text{Eq.~}\eqref{#1}}}

\newcommand{\secref}[1]{Sec.~\ref{#1}}
\newcommand{\tabref}[1]{Tab.~\ref{#1}}

\usepackage{amsmath, amssymb, amsfonts, epsfig, xspace}

\newcommand{\norm}[1]{\left\lVert#1\right\rVert}
\newcommand{\diff}[1]{\ensuremath{\operatorname{d}\!{#1}}}
\newcommand{\st}[0]{\enskip \text{s.t.:} \enskip}
\newcommand{\abs}[1]{\left\lvert#1\right\rvert}    
    
\newcommand{\prox}[0]{\text{prox}}

\newcommand{\mps}[0]{\ensuremath{\frac{\text{m}}{\text{s}} \;}} 
\newcommand{\symspd}[1]{\ensuremath{\mathcal{P}(#1)}} 
\newcommand{\localApproxOp}[0]{\ensuremath{U}} 
\DeclareMathOperator*{\softmin}{\text{smin}}

\newcommand{\caseif}[0]{\; \text{if} \;} 
\newcommand{\caseelse}[0]{\; \text{else} \;} 
 
\makeatletter
\newcommand{\@chapapp}{\relax}%
\makeatother

\usepackage[title,toc,titletoc,header]{appendix}

\makeatletter
\def\thanks#1{\protected@xdef\@thanks{\@thanks
        \protect\footnotetext{#1}}}
\makeatother

\newcommand{\TheMethod}{\textsc{PIEMAP}}
\title{\TheMethod: Personalized Inverse Eikonal
Model from cardiac Electro-Anatomical Maps 
\thanks{
	This research was supported by the grants F3210-N18 and I2760-B30
	from the Austrian Science Fund (FWF) and BioTechMed Graz flagship award 
	"ILearnHeart", as well as ERC Starting grant HOMOVIS, No. 640156. This work was also financially supported by the Theo Rossi di Montelera Foundation, the Metis Foundation Sergio Mantegazza, the Fidinam Foundation, the Horten Foundation and the CSCS—Swiss National Supercomputing Centre production grant s778.
}
}
\author{Thomas Grandits\inst{1, 4}, Simone Pezzuto\inst{2}, Jolijn M. Lubrecht\inst{2}, Thomas Pock\inst{1, 4}, Gernot Plank\inst{3, 4}, Rolf Krause\inst{2}}
\institute{Institute of Computer Graphics and Vision\\
Graz University of Technology\\
\email{\{thomas.grandits,pock\}@icg.tugraz.at}
\and
Center for Computational Medicine in Cardiology,\\
Institute of Computational Science,\\
Universit\`a della Svizzera italiana,\\
Lugano, Switzerland\\
\email{\{simone.pezzuto,jolijn.marieke.lubrecht,rolf.krause\}@usi.ch}
\and 
Institute of Biophysics\\
Medical University of Graz\\
\email{gernot.plank@medunigraz.at}
\and
BioTechMed-Graz, Austria
}

\date{\today}

\begin{document}

\maketitle

\begin{abstract}
    Electroanatomical mapping, a keystone
    diagnostic tool in cardiac electrophysiology studies,
    can provide high-density maps of the local electric
    properties of the tissue. It is therefore
    tempting to use such data to better individualize current
    patient-specific models of the heart
    through a data assimilation procedure
    and to extract potentially insightful information
    such as conduction properties.
    Parameter identification for state-of-the-art cardiac models is
    however a challenging task.
    
    In this work, we introduce a novel inverse problem for inferring
    the anisotropic structure of the conductivity tensor,
    that is fiber orientation and conduction velocity along and across fibers,
    of an eikonal model
    for cardiac activation.
    %
    The proposed method, named \TheMethod, performed robustly with synthetic data and showed promising results
    with clinical data.
    These results suggest that \TheMethod\ could be a useful supplement in future clinical workflows
    of personalized therapies.
\end{abstract}

\section{Introduction}
Patient-specific modeling in cardiac electrophysiology has nowadays reached the status of a clinically feasible tool for assisting the cardiologist during the therapeutic intervention.  As models became more mature, and thanks to the increasingly availability of high-resolution data such as high-density electroanatomic maps (EAMs), parameter identification has emerged as a key topic in the field.

A high-density EAM is composed by a large number of contact recordings (1000 points or more), each with local electrogram and spatial information. Activation and conduction velocity maps, for instance, can be derived by combining electric and geometric data.  Conductivity parameters in a propagation model may therefore be adapted to reproduce such maps for model personalization.
%
%

The reconstruction of conduction velocity maps is generally based on local approaches
\cite{cantwell_techniques_2015,Coveney2020}.
In these methods, the local front velocity is estimated from an appropriate interpolation of the local activation time (LAT). Anisotropic conductivity can be deduced from front velocity and prior knowledge on fiber structure (rule-based or atlas-based), or by combining multiple activation maps~\cite{roney_technique_2019}.

Although being computationally cheap, these local methods may miss effective mechanisms for global consistency in electric wave propagation models, and may introduce artefacts in conduction velocity due to front collisions or breakthroughs.
A different approach, also adopted in this work, relies instead on the
(possibly strong) assumption that a calibrated model for the cardiac
activation can reproduce the measured activation with sufficient accuracy.
The electric conductivity in the model is eventually identified through an optimization procedure
aiming at minimizing the mismatch between the model output and the
collected data. The model can either be enforced pointwise, yielding for instance PDE-constrained optimization~\cite{barone2020}, or act as a penalization term~\cite{sahli_costabal_physics-informed_2020}.

To the best of our knowledge, however, the problem of estimating \emph{simultaneously} distributed fiber architecture and conduction velocities from sparse contact recordings has never been attempted before with either approaches. 
In this work we aim to bridge this gap by proposing a novel method to extract from a single EAM the full electric conductivity
tensor, with the only assumption of symmetric positive-definiteness (s.p.d.) of the tensor field. Local fiber
orientation and conduction velocities are then deduced
from the eigendecomposition of the conductivity tensor.
As forward model, we consider
the anisotropic eikonal model, which is a good compromise between physiological accuracy
and computational cost~\cite{franzone_spreading_1993}.
The corresponding inverse model, employing Huber 
regularization, a smooth total variation approximation, to stabilize the reconstruction
and log-Euclidean metric in the parameter space to ensure
s.p.d.\ of the tensor field,
is solved by an iterative quadratic approximation strategy 
combined using
a Primal-Dual optimization algorithm.
Finally, we extensively test the algorithm with synthetic
and clinical data for the activation of the atria, represented as a 2-D manifold, showing promising results for its clinical applicability.

%
%
\section{Methods}
\subsection{Forward Problem}
The anisotropic eikonal equation describes the activation times $u$ of a wave propagating with direction-dependent velocity. Given a smooth 2-D manifold $\Omega\subset\mathbb{R}^3$, the equation reads as follows
\begin{equation}
\begin{split}
    \sqrt{\left<D(\mathbf{x}) \nabla\!_\mathcal{M} u(\mathbf{x}), \nabla\!_\mathcal{M} u(\mathbf{x})\right>} = 1 \st &\forall \mathbf{x} \in \Omega: D_3(\mathbf{x}) \in \symspd{3} \\ &\forall \mathbf{x} \in \Gamma_0\subset\partial\Omega: u(\mathbf{x}) = 0,
\end{split}
    \label{eq:aniso_eikonal}
\end{equation}
with $\Gamma_0$ representing the domain of fixed activation times, $\symspd{n}$ being the space of $n\times n$
symmetric positive definite matrices and $\nabla\!_\mathcal{M}$ being the surface gradient. The conductivity tensor $D(\mathbf{x})$ specifies the conduction velocity in the propagation direction, that is $\sqrt{\left<D(\mathbf{x}) \mathbf{k}, \mathbf{k}\right>}$ is the velocity of the wave at $\mathbf{x}\in\Omega$ in direction $\mathbf{k}$, $\mathbf{k}$ unit vector.

To solve \eqreffull{eq:aniso_eikonal}, we based our algorithm---purely implemented in TensorFlow to allow for automatic gradient computation through back-propagation---on the Fast Iterative Method (FIM) 
for triangulated surfaces~\cite{fu_fast_2011}. The only fixed assumed point in \TheMethod\ $\{\mathbf{x}_0\} = \Gamma_0$ is the chosen earliest activation site, which is assumed to be the earliest of all measured activation points.
We use a slightly altered fixed-point iteration $\mathbf{u}^{k+1} = F(\mathbf{u}^k)$ 
which iteratively updates the activation times. For sake of simplicity, we keep denoting by $u(\mathbf{x})$ the piecewise linear interpolant of the nodal values $u_i$ and by $D$ a piecewise constant tensor field on the triangulated surface.  The approximated solution of the equation~\eqref{eq:aniso_eikonal}, henceforth denoted by $\operatorname{FIM}_D(\mathbf{x})$ is then the unique fixed-point of the map $F$.
Specifically, the map $F$ updates each of the nodal values as follows: 
\begin{equation}
u^{k+1}_i =
\left\{\begin{aligned}
&0, & \mbox{if $\mathbf{x}_i\in\Gamma_0$},\\
& \softmin^\kappa_{T_j \in \omega_i}
 \softmin^\kappa_{\mathbf{y}\in e_{i,j}}
 \biggl\{
 u(\mathbf{y})
 + \sqrt{\bigr< D^{-1}_j
 (\mathbf{y} - \mathbf{x}_i),
 \mathbf{y} - \mathbf{x}_i\bigr>}
 \biggr\},
& \mbox{otherwise}.
\end{aligned}\right.
    \label{eq:fim_update}
\end{equation}
%
%
where $\omega_i$ is the patch of triangles $T_j$ connected to the vertex $\mathbf{x}_i$, $e_{i,j}$ is the edge of the triangle $T_j$ opposite to the vertex $\mathbf{x}_i$, $D_j = D|_{T_j}$ and $\softmin^k$ being the soft-minimum function, defined as
$\softmin^\kappa (\mathbf{x}) = -\frac{1}{\kappa} \log ( \sum_i \exp ( -\kappa x_i ) )$. 

Differently from the classic FIM method~\cite{fu_fast_2011}, we concurrently update \emph{all} the nodes, that is the map $F$ is applied in parallel to each node and not just on a small portion of ``active'' nodes. We also replaced the $\min$-function of the original FIM-algorithm by the 
soft-minimum $\softmin^\kappa$ to ensure a limited degree of smoothness of the function
and avoid discontinuities in the gradient computation. 

\subsection{\TheMethod\ and inverse problem}
\label{sec:min_problem}

\TheMethod\ implements an inverse problem in which the optimal conductivity tensor field $D$ for \eqreffull{eq:aniso_eikonal} is selected such that the mismatch between recorded activation times $\hat{u}(\mathbf{x})$ and the simulated activation times $\operatorname{FIM}_D(\mathbf{x})$ on the measurement domain $\Gamma\subset\Omega$
is minimized in the least-squares sense.

\begin{figure}[h]
	\centering
		\includegraphics[width=.8 \textwidth, trim={1cm, 0cm, 1cm, 0cm}, clip]{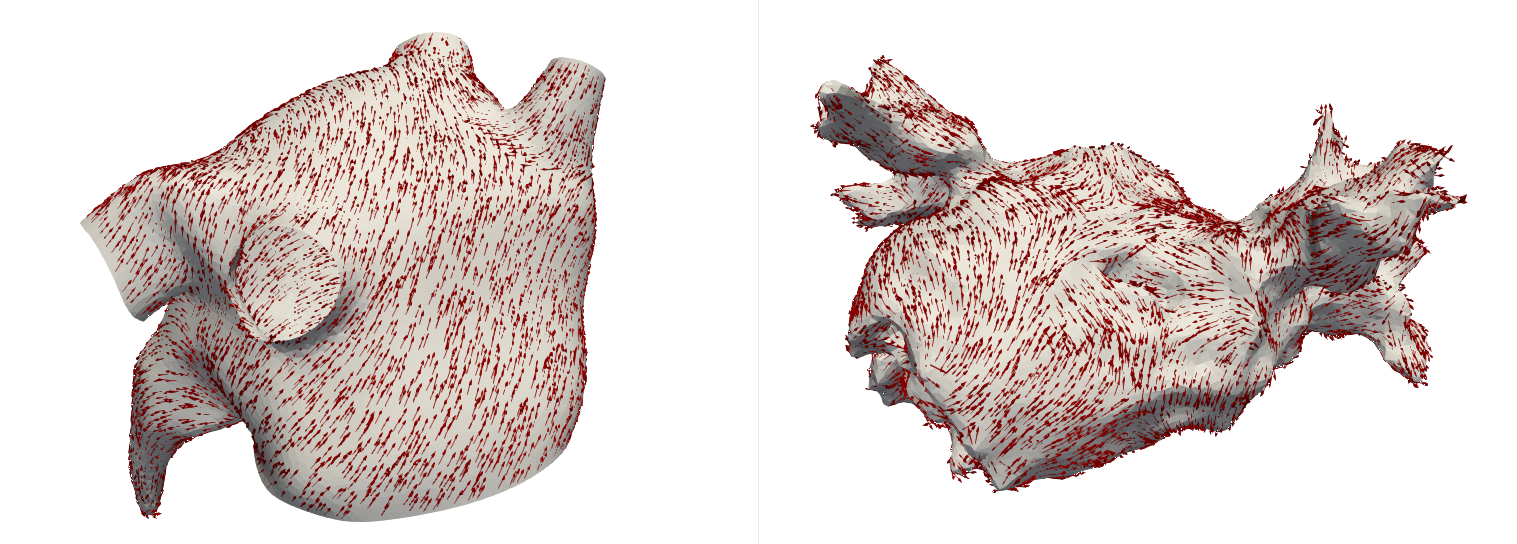}
	%
	\caption{Generated local bases on the atria manifold models.}
	\label{fig:local_bases}
\end{figure}

In principle, after accounting for symmetry, $D(x)$ has 6 component to be identified for every $x\in\Omega$.  Since the dynamic of the wave propagation is bound to the 2-D manifold, however, the component normal to the surface does not influence the solution. We therefore define $D(\mathbf{x})$ as follows:
\begin{equation}
D(\mathbf{x}) = P(\mathbf{x}) \begin{pmatrix}
\tilde{D}(\mathbf{x}) & \mathbf{0} \\
\mathbf{0}^\intercal & 1
\end{pmatrix}
P(\mathbf{x})^\intercal,
\label{eq:conductivity3D}
\end{equation}
where $\tilde{D}(\mathbf{x})\in\symspd{2}$ and $P(x)$ is a rotation from the canonical base in $\mathbb{R}^3$ to a local base $\{ \mathbf{v}_1(\mathbf{x}), \mathbf{v}_2(\mathbf{x}), \mathbf{n}(\mathbf{x}) \}$. The local base at $\mathbf{x}\in\Omega$ is such that $\mathbf{n}(\mathbf{x})$ is the normal vector to the surface and $\mathbf{v}_1(\mathbf{x})$, $\mathbf{v}_2(\mathbf{x})$ span the tangent space.  In such way, the dimension of the parameter space is reduced from 6 to only 3.
Any basis $\mathbf{v}_1(\mathbf{x})$ in the tangent space is valid, but we compute a smooth basis by minimizing the variation across the manifold:
\begin{equation*}
    \min_{\mathbf{v}_1} \int_{\Omega} \norm{\nabla \mathbf{v}_1(\mathbf{x})}_2^2 \diff{\mathbf{x}} \st \forall \mathbf{x} \in \Omega: \norm{\mathbf{v}_1(\mathbf{x})} = 1,
\end{equation*}
to ensure a meaningful result through the later introduced regularization term. The computed local bases, used in all experiments throughout this paper, are shown in \figref{fig:local_bases}.
%
Finally, we consider the Log-Euclidean metric \cite{arsigny_geometric_2007} for ensuring a s.p.d.\ tensor field: given $\mathbf{d}(\mathbf{x}) \in\mathbb{R}^3$, $\mathbf{x}\in\Omega$, we set $\tilde{D}$ as follows:
\begin{equation}
\tilde{D}(\mathbf{x}) :=
\exp \begin{pmatrix}
d_1(\mathbf{x}) & d_2(\mathbf{x}) \\
d_2(\mathbf{x}) & d_3(\mathbf{x})
\end{pmatrix},
\label{eq:logmetric}
\end{equation}
where the matrix exponential is computed from the eigendecomposition.
In particular, the admissible set $\mathbb{R}^3$ is mapped through \eqref{eq:conductivity3D} and \eqref{eq:logmetric} to $\symspd{3}$.

%
%

The inverse problem, therefore, consists in finding the vector field $\mathbf{d} \in \mathbb{R}^3$, 
which minimizes the following objective function:
\begin{equation}
    \min_{\mathbf{d}}
    \underbrace{\frac{1}{2} \int_{\Gamma} \left( \operatorname{FIM}_{D(\mathbf{d})}(\mathbf{x}) - \hat{u}(\mathbf{x}) \right)^2 \diff{\mathbf{x}}}_{
    \localApproxOp(\mathbf{d})} +
    \underbrace{\lambda \int_{\Omega} H_{\varepsilon}( \nabla\!_{\mathcal{M}} \mathbf{d})  \diff{\mathbf{x}}}_{\operatorname{TV}_{\varepsilon,\lambda}(\mathbf{d})},
    \label{eq:orig_problem}
\end{equation}
where $\operatorname{TV}_{\varepsilon,\lambda}(\mathbf{d})$
is a smooth total variation (TV) regularization term which alleviates the ill-posedness of the problem.
Specifically, $H_\varepsilon$ is the Huber function:
\begin{equation}
    H_{\varepsilon}(\mathbf{x}) = 
    \begin{cases}
     \frac{1}{2 \epsilon} \abs{\mathbf{x}}^2, & \caseif \abs{\mathbf{x}} \le \varepsilon, \\
     \abs{\mathbf{x}} - \frac{1}{2} \varepsilon, & \caseelse
    \end{cases}
\label{eq:huber}
\end{equation}

We set $\epsilon = \num{5e-02}$ for our experiments, while the optimal choice of regularization parameter $\lambda$ is obtained by using a cross-validation approach.
\subsection{Forward-Backward Splitting and Numerical Solution}
The computational complexity of solving \eqreffull{eq:orig_problem} is dominated by the time for 
computing $\text{FIM}_D$ and $\nabla_\mathbf{d} \text{FIM}_D$. The implementation of $\text{FIM}_D$ in TensorFlow allows for an efficient
computation of $\nabla_\mathbf{d} \text{FIM}_D$ via backpropagation on a graphical processing unit (GPU).
While the minimization of the residual $U(\mathbf{d})$ is usually achieved very quickly, at least when $\hat{u}$ is a (possibly corrupted) solution of Eq.~\eqref{eq:aniso_eikonal},
the TV term tends to increase the number of needed iterations for convergence. In order to increase the convergence rate, we apply the principle of the Fast
Iterative Shrinking and Thresholding Algorithm (FISTA)~\cite{beck_fast_2009},
quadratically bounding the non-linear, non-convex 
function $\localApproxOp$ around the current point $\mathbf{d}_k$:
\begin{equation}
    \localApproxOp(\mathbf{d}) \le \localApproxOp(\mathbf{d}_k) + \left<\nabla_\mathbf{d} \localApproxOp(\mathbf{d}_k), \left( \mathbf{d} - \mathbf{d}_k \right) \right> + \frac{L}{2} \norm{\mathbf{d} - \mathbf{d}_k}_2^2 =: G(\mathbf{d}).
    \label{eq:descent_lemma}
\end{equation}
The bounding function $G(\mathbf{d})$ is convex, hence has a unique minimum $\bar{\mathbf{d}} = \mathbf{d}_k - L^{-1} \nabla_\mathbf{d}\localApproxOp(\mathbf{d}_k)$. As
$\text{TV}_{\varepsilon,\lambda}(\mathbf{d})$ is also convex, we obtain the following convex
minimization problem:
\begin{equation*}
    \min_{\mathbf{d} } \frac{L}{2} \Vert \mathbf{d} - \bar{\mathbf{d}}\Vert_2^2 + \text{TV}_{\varepsilon, \lambda}(\mathbf{d}).
\end{equation*}
Iteratively solving this class of problems along with an acceleration term is usually referred to as FISTA. 
We recast the problem into a convex-concave saddle-point problem: 
\begin{equation}
    \min_{\mathbf{d} } \max_{\mathbf{p}} \frac{L}{2} \norm{\mathbf{d} - \bar{\mathbf{d}}}_2^2 + \left<\nabla_{\mathcal{M}} \mathbf{d}, \mathbf{p} \right> - \text{TV}_{\varepsilon, \lambda}^* (\mathbf{p})
    \label{eq:local_pd}
\end{equation}
which can be solved using the Primal-Dual algorithm \cite{chambolle_first-order_2011} given by:
\begin{equation}
\begin{cases}
 &\mathbf{d}^{i+1} = \prox_{\tau G} (\mathbf{d}^i - \tau \nabla^{*}_{\mathcal{M}} \mathbf{d}^i) \\
 &\mathbf{d}_{\Theta} = \mathbf{d}^{i+1} + \theta \left( \mathbf{d}^{i+1} - \mathbf{d}^i \right) \\
 &\mathbf{p}^{i+1} = \prox_{\sigma \text{TV}_{H_{\epsilon}, \mathcal{M}}^*} \left( \mathbf{p}^k + \sigma \nabla_{\mathcal{M}} \mathbf{d}_{\Theta} \right)
 \end{cases}
\end{equation}
\todo{Check divergence sign}
with 
\begin{equation*}
\begin{split}
\hat{\mathbf{d}} = \prox_{\tau U} (\tilde{\mathbf{d}}) &= \left( \tilde{\mathbf{d}} + \tau L \bar{\mathbf{d}} \right) / \left( \tau L + 1 \right) \\
\hat{\mathbf{p}} = \prox_{\sigma \text{TV}_{\epsilon, \mathcal{M}}^*} (\mathbf{p}) &\Leftrightarrow 
\hat{\mathbf{p}}_j =
\begin{cases}
 \frac{\bar{\mathbf{p}}_j}{\abs{\bar{\mathbf{p}}_j} / \lambda} &\caseif \abs{\bar{\mathbf{p}}_j} > 1 \\
 \bar{\mathbf{p}}_j &\caseelse
\end{cases}
\end{split}
\end{equation*}
for $\bar{\mathbf{p}}_j = \frac{\mathbf{p}_j}{\sigma \epsilon / \lambda + 1}$, $\theta = 1$ 
and $\tau \sigma \norm{\nabla_{\mathcal{M}}}_2^2 \le 1$.
%
The parameter $L$ in \eqreffull{eq:descent_lemma}, usually challenging to evaluate, is computed through 
a Lipschitz backtracking algorithm~\cite{beck_fast_2009}.
\section{Experiments}
\label{sec:experiments}
For the evaluation of \TheMethod, we first assessed its effectiveness on reconstructing known conduction velocity and fibers on a realistic human left atrium (LA) model, also in the presence of white noise and heterogeneity. 
The LA model was generated from MRI data of a patient, with the fibers semi-automatically assigned as described previously~\cite{gharaviri20}. Fiber and transverse velocity were set to $0.6$ \mps and $0.4$ \mps respectively for the entire LA, except for the low conducting region, where we used $0.2$ \mps for both fiber and transverse velocity.
We tested \TheMethod\ both in the case of fully anisotropic and in the case of isotropic conduction. In the latter case, in particular, we compared \TheMethod\ to existing methods for the evaluation of conduction velocity, namely a local method~\cite{cantwell_techniques_2015} and EikonalNet~\cite{sahli_costabal_physics-informed_2020}, a Physics Informed Neural Network (PINN) method.
In a second set of experiments, we eventually applied \TheMethod\ to clinically acquired data, in the form of high-density EAM.

All experiments were run on a desktop machine with an Intel Core i7-5820K CPU with 
6 cores of each 3.30GHz, 32GB of working memory and a NVidia RTX 2080 GPU.
All examples were optimized
for 2000 iterations, with each iteration taking about $1.8$ seconds, totalling into 
a run-time of approximately 1 hour for one optimization.

\subsection{Numerical assessment}
\label{sec:experiments_la}


All the experiments were performed on a human, cardiac magnetic resonance (CMR)-derived left atrium model, with semi-automatically placed fiber directions based on histological studies.
The ground-truth (GT) solution was computed with a single earliest activation site using \eqreffull{eq:aniso_eikonal}, and with a low-conducting area being close to the left atrial appendage.  Different levels of independent and identically distributed (i.i.d.) Gaussian noise with standard deviation $\sigma_{\mathcal{N}}$ were tested. The measurement domain was a set of 884 points uniformly distributed across the atrium.
The reconstruction root-mean-square error (RMSE) with respect to GT was evaluated in terms of conduction velocity (m/s), propagation direction and, only for \TheMethod, fiber-angle error.
To evaluate the results, we compute the front direction and fiber direction unit vectors, denoted
as $\mathbf{e}$ and $\mathbf{f}$ respectively.
The front and fiber angle-errors are then defined as 
$\alpha_{\mathbf{e}} = \arccos \left<\mathbf{e}, \mathbf{e}_{\text{GT}}
\right>  \in \left[0, 180^{\circ} \right)$ and
$\alpha_{\mathbf{f}} \arccos \abs{\left<\mathbf{f}, \mathbf{f}_{\text{GT}} 
\right>}  \in \left[0, 90^{\circ} \right)$. 
The velocity errors in propagation direction are then  
$v \cos \left( \alpha_{\mathbf{x}} \right) - v_{\text{GT}}$
for computed velocity $v$ and exact velocity $v_{\text{GT}}$, both in the front and fiber direction.

\def\arraystretch{1.2}
\setlength{\tabcolsep}{.4em}
\begin{table}[htb]
\centering
\caption{Comparison of the front-velocity/front-angle error of \TheMethod\ with the local approach \cite{cantwell_techniques_2015} and EikonalNet~\cite{sahli_costabal_physics-informed_2020}, assuming different noise levels
for the in-silico models. Errors in $\frac{m}{s}$/degree. 
In the last column, we compare the fiber velocity error and fiber angle error.
}
\begin{tabular}{cc|c|c|c||c|}
\cline{3-6}
     &   & \multicolumn{3}{c||}{\textbf{Error in Propagation Direction}} & \multicolumn{1}{c|}{\textbf{Fiber Error}} \\ \cline{3-6} 
     &   & \TheMethod  & Local Method & EikonalNet & \TheMethod \\ \hline
     \multicolumn{1}{|l|}{\multirow{4}{*}{\rotatebox[origin=c]{90}{$\sigma_{\mathcal{N}}$/PSNR}}} & 0ms/$\infty$ dB & \textbf{\num{1.9926e-01}}/\num{1.0580e+01} & \textbf{\num{1.9901e-01}}/\num{2.2950e+01} & \num{5.29e-01}/\textbf{\num{9.20e+00}} & \num{2.4804e-01}/\num{3.8343e+01}\\ \cline{2-6}
\multicolumn{1}{|l|}{}       & 0.1ms/64.1 dB & \textbf{\num{1.9359e-01}}/\num{1.0606e+01} & \num{2.0228e-01}/\num{2.3170e+01} & \num{3.96e-01}/\textbf{\num{9.10e+00}} & \num{2.4805e-01}/\num{3.8457e+01}\\ \cline{2-6}
 \multicolumn{1}{|l|}{}       & 1ms/43.9 dB & \textbf{\num{1.9615e-01}}/\textbf{\num{1.1033e+01}} & \num{2.1177e-01}/\num{2.3566e+01} & \num{4.88e-01}/\num{1.46e+01}& \num{2.4877e-01}/\num{3.8594e+01}\\ \cline{2-6}
\multicolumn{1}{|l|}{}       & 5ms/29.9 dB & \textbf{\num{2.5130e-01}}/\textbf{\num{1.9936e+01}} & \num{2.9125e-01}/\num{3.0197e+01} & \num{1.24e+00}/\num{4.94e+01} & \num{2.5884e-01}/\num{4.0136e+01}\\ \hline
\end{tabular}
	\label{tab:comparison_tab}
\end{table}

Results are reported in~\tabref{tab:comparison_tab}.
All methods correctly captured the low conduction region.
\TheMethod\ compared favourably to the local method at all noise levels in terms of absolute conduction velocity. EikonalNet shows a slightly more accurate front angle error, which is counteracted by the considerably high front velocity error, both compared to our and the local method. 

Overall, \TheMethod\ had the benefit over EikonalNet that the GT was generated with the anisotropic eikonal model, and thus it is in theory possible to reproduce the data exactly with a zero noise level. In the local method no model assumption is made. Interestingly, the error in front direction for the local method could be linked to the fact that, in the presence of anisotropic conduction, propagation direction and $\nabla\!_\mathcal{M} u$ differ. For instance, a circular propagation from the source $x_0$ satisfies Eq.~\eqref{eq:aniso_eikonal} with $u(x)=\sqrt{\left<D^{-1}(x-x_0),(x-x_0)\right>}$, thus $\nabla\!_\mathcal{M} u$ differs from $x-x_0$, which is the propagation direction.  In the local method, $\nabla\!_\mathcal{M} u$ is used to establish such direction.  In EikonalNet, results were less robust to noise. A plausible explanation is that training Neural Networks does not always yield the same results, as multiple local minima might be present. Therefore, error can be slightly lower or higher depending on the initial conditions.
In terms of computational time, \TheMethod\ was comparable to EikonalNet, but significantly slower than the local method.
\begin{figure}[htb]
    \centering
    \includegraphics[width=.9\textwidth, 
     trim={0.5cm 4.5cm 2.5cm 5.5cm}, 
    clip]{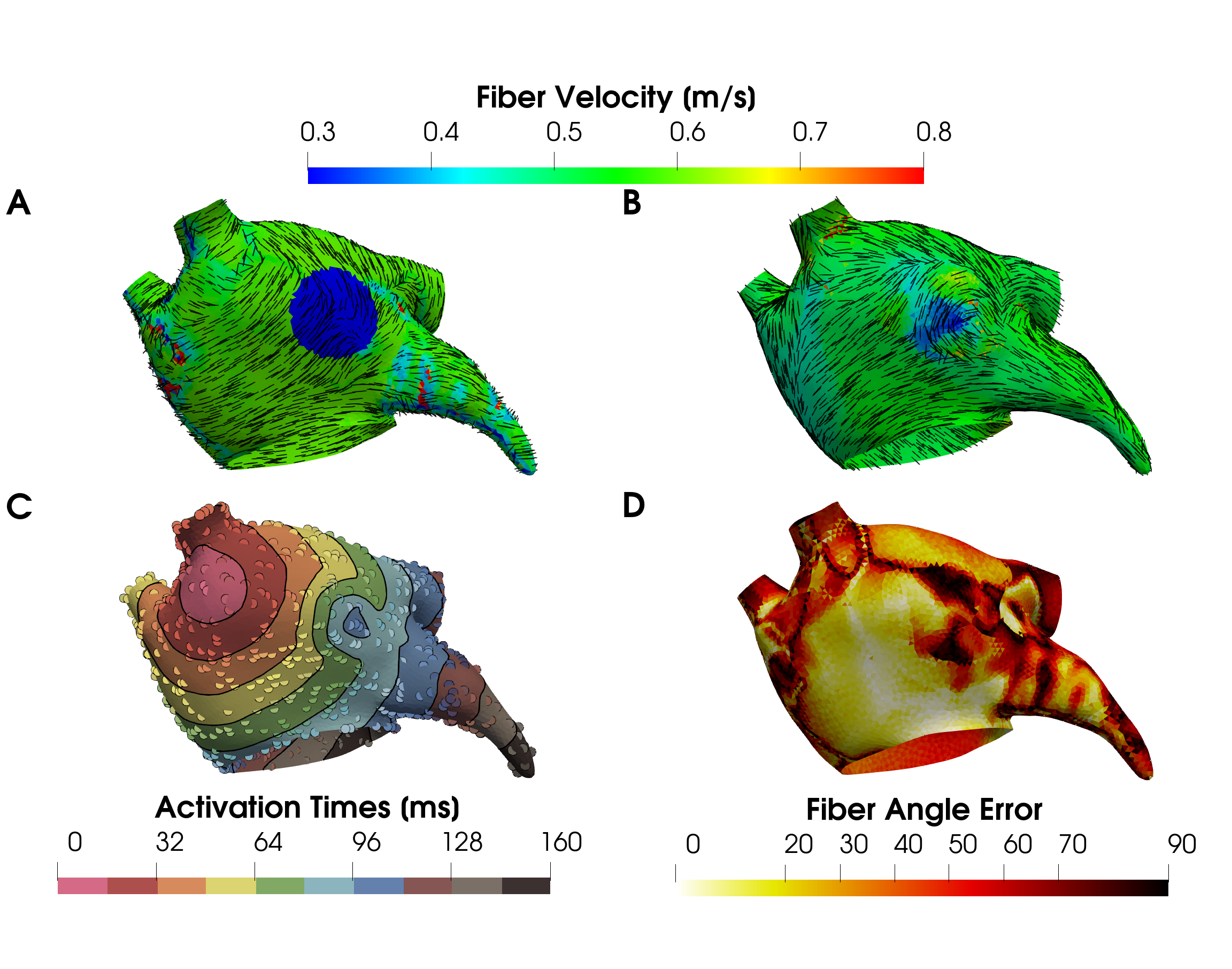}
    \caption{Results of our method on the noise-less LA-model with known ground-truth fiber orientation and velocity (A). The scarred region is correctly activated at a later time by a combination of reducing fiber-velocity, as well as aligning fibers along the contour lines (B/C). The activation map can faithfully capture the observed measurement points $u(\mathbf{x} \in \Gamma)$, marked as colored dots (C) and matches the GT-model's activation closely (not shown). Fiber alignment of our model (B) shows mostly errors around the mitrial valves, the scarred region and pulmonary veins (D), as well as regions of high curvature. 
    Best viewed online.}
    \label{fig:fibers_LAmodel}
\end{figure}

Regarding the reconstruction of fiber directions (see \figref{fig:fibers_LAmodel}), we observed a very good performance for the fiber and cross-fiber velocity, and a reasonable reconstruction for the direction. In particular, reconstruction in fiber direction was poor around the boundaries (mitral ring and pulmonary veins, where fibers run parallel to the opening) and in the scarred region, which attribute the most to the fiber angle error in \tabref{tab:comparison_tab}. The distribution of fiber angle errors is a slightly left-skewed uniform distribution (not shown), indicating that the chosen smooth basis along with a simple TV prior can provide resonable results with respect to the activation timings. Still, it may not be sufficient to account for the partly complicated fiber orientation, especially in areas of high-curvature of the mesh or sudden changes of fiber orientation on the endocardium as an effect of the volumetric structure of the atria, such as is the case for the mitrial valve. Physiological priors will need to be considered in the future for this purpose.

\subsection{Application to real clinical data}
\label{sec:acquisition_eam}

In a patient candidate to ablation therapy, a high-density activation map along with a 3D patient-specific atrial model was acquired with an EAM system 
(Catheter: Pentaray\textregistered{} System: CARTO\textregistered{} 3 System, Biosense Webster). 
The recordings encompassed roughly 850 ``beats'' of 2.5 sec including both the electrode position in 3D space and the unipolar electrogram (1\,kHz). Recordings that were deemed to be untrustworthy due to 1) insufficient contact, 2) sliding of the electrode in 3D space \textgreater 1 cm, 3) correspondence to a inconsistent surface P-wave, 4) minimal unipolar amplitude, were excluded automatically from the study. To avoid degenerated triangles with acute angles, sometimes created by the EAM recordings, we used PyMesh\footnote{\url{https://github.com/PyMesh/PyMesh}} to postprocess the mesh.
A further
manual pre-processing of the signals was eventually performed for a correct detection of the local activation time (steepest negative deflection in the unipolar signal) in the last beat and compared to local bipolar signals for confirmation. Distribution of points was uneven across the LA, as many points were located around the pulmonary veins (PVs). 

\begin{figure}[tbh]
	\centering
	\includegraphics[width=.9\textwidth, trim={0.5cm 3.5cm 0.5cm 3.5cm}, clip]{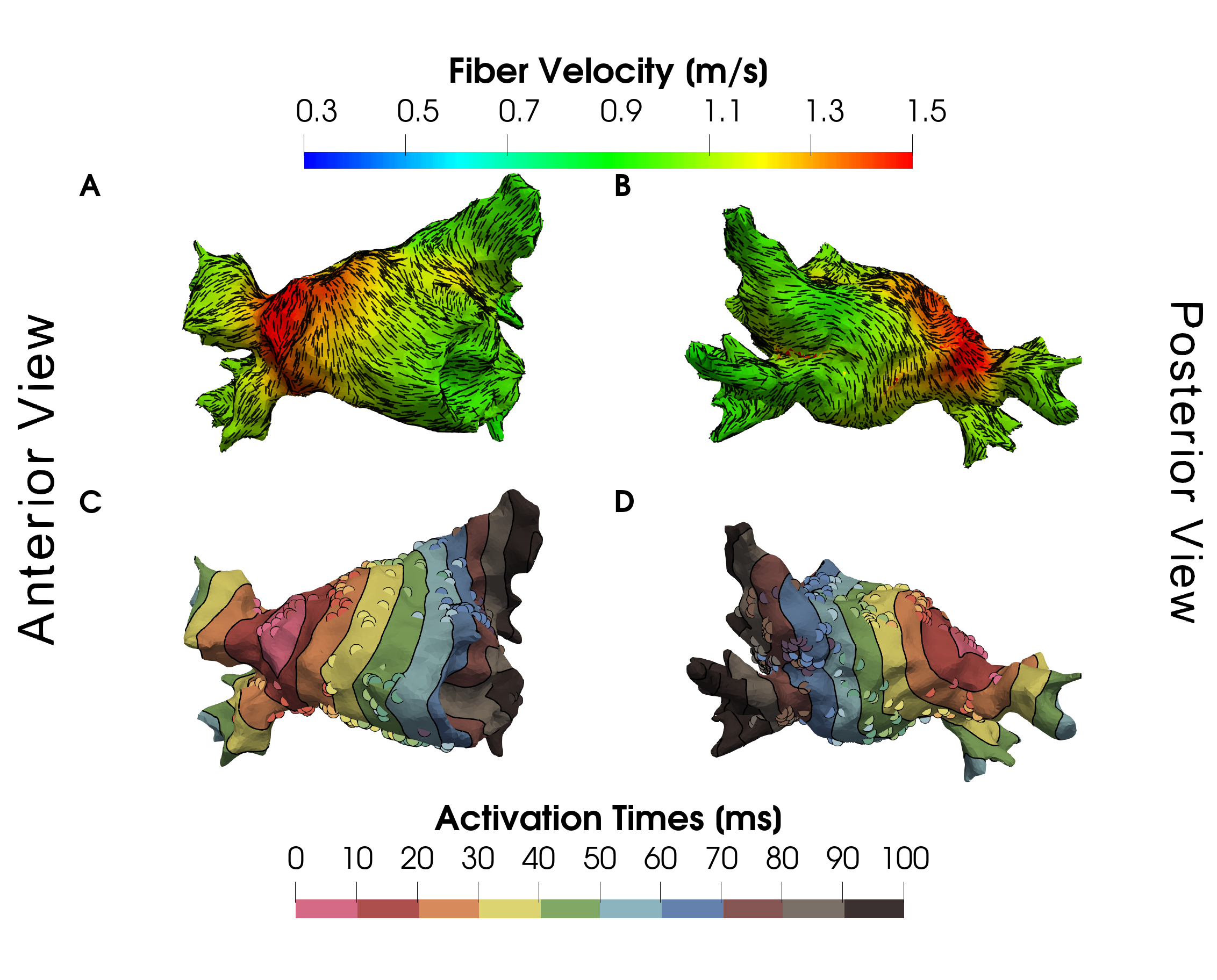}
	\caption{Anterior (A/C) and posterior view (B/D) of \TheMethod's results on a patient's left atria. The panels A and B show the found fiber direction and fiber velocity, while
	the panels C and D show the activation map along with the actual measured points on top, similar to \figref{fig:fibers_LAmodel}.}
	\label{fig:realdata}
\end{figure}

Of the remaining valid 565 beats, randomly chosen $80\%$ ($452$ points) were used to optimize \eqreffull{eq:orig_problem}, while the remaining $20\%$ were used as a cross-validation set to find the optimal regularization parameter $\lambda$. \figref{fig:realdata} shows the fiber velocity, orientation and activation map after the optimization. The cross-validation error over several values of $\lambda$ is shown in \figref{fig:EAM_lambda_study}, which lead us to the used value of $\lambda$. The best cross-validation error lead to a relatively smooth fiber velocity field, with velocities ranging up to $1.5 \mps$ in the initiation region, probably a consequence of choosing only one mesh node as an initiation site when in reality the initiation site is larger or composed of multiple sites. A speed-up of propagation near the atrial wall can often be witnessed and is compensated in our model by an overall higher fiber-velocity.

\sisetup{scientific-notation=true, exponent-product = \cdot, round-precision=1}

\bgroup
\sisetup{scientific-notation=true, exponent-product = \cdot}
\begin{figure}[h]
    \centering
    \subfloat{
      \includegraphics[width=.8\textwidth, trim={0cm 1cm 0cm 0.0cm}, clip]{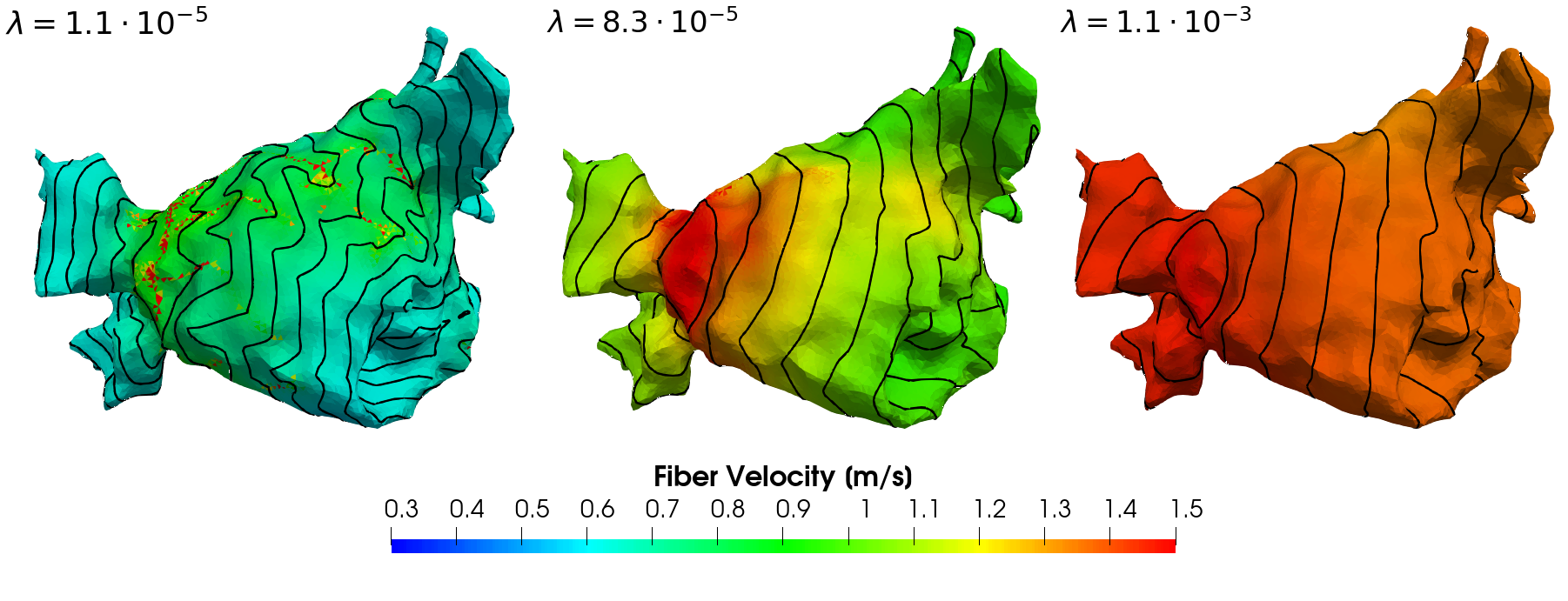}
    }
    
    \subfloat{
      \includegraphics[width=.7\textwidth, trim={1.5cm 0.0cm 1.5cm 0.5cm}, clip]{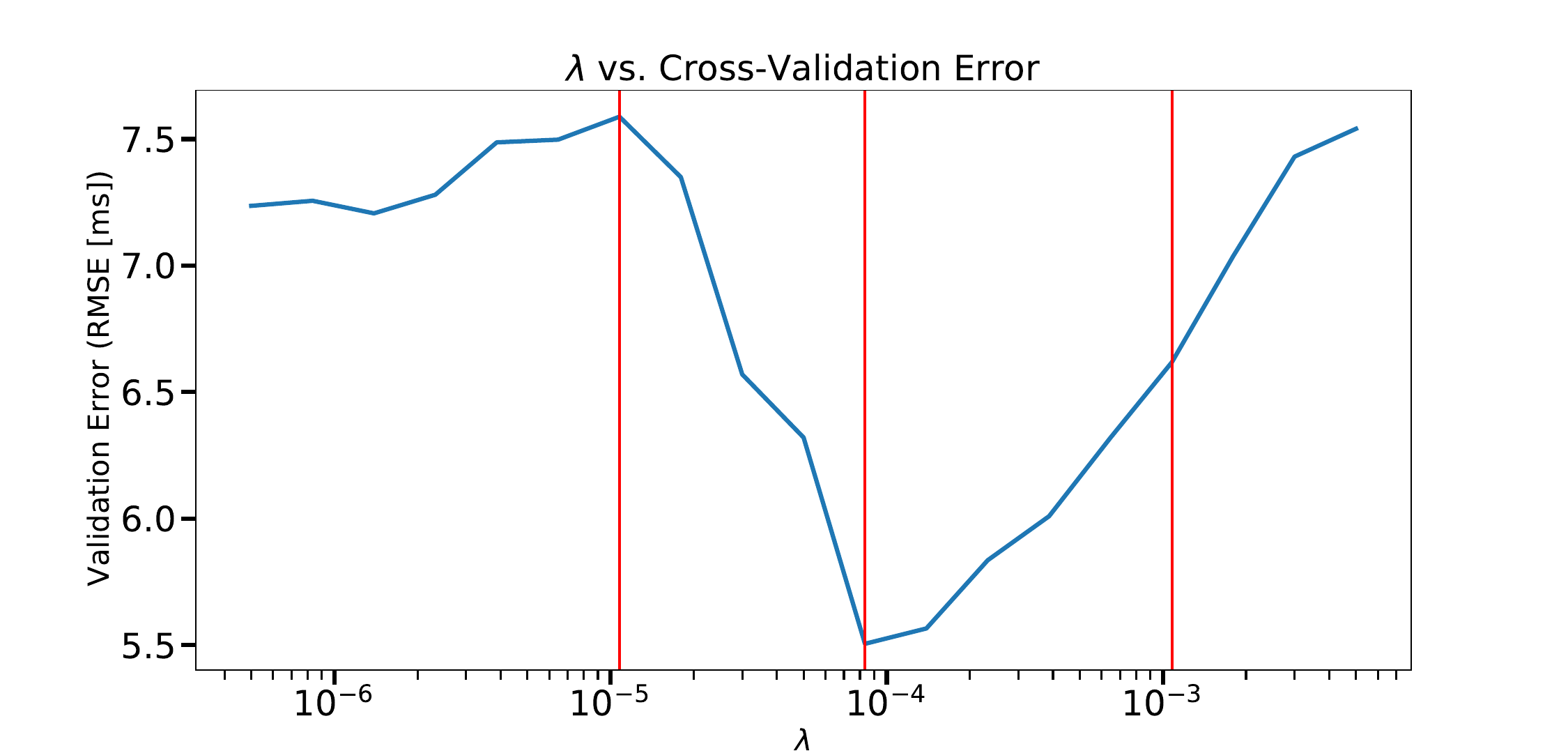}
    }
    \caption{Influence on the final cross-validation error when varying over $\lambda$ for the optimization on the EAM recordings. 
    The black contours are the isochrones of the modelled activation. Shown are the results of
    three different $\lambda$ values: 
    The left result is the least regularized with $\lambda = \num{1.08e-5} $, while the result on the right side is heavily regularized with $\lambda = \num{1.08e-03}$. A compromise is the figure in the middle with $\lambda = \num{8.34e-5}$, but
    finding the physiologically most plausible  value for $\lambda$ is a non-trivial task since
    we do not know the true distribution of velocities and activations in the atria.
    }
    \label{fig:EAM_lambda_study}
\end{figure}
\egroup

\section{Discussion \& Conclusion}
In this paper, we proposed \TheMethod, a global method to reconstruct the conductivity tensor (fiber direction, fiber- and cross-fiber velocity) of an anisotropic eikonal model from sparse measurement of the activation times.  We compared our method to existing approaches for determining the conduction velocity map from the same data (a local method and a PINN method) and we demonstrated its effectiveness in a real application.

Our method showed promising results on atrial electrical data, acquired using an EAM system, but may be used with any electrical measurements, mapped to a manifold. In \secref{sec:experiments_la}, we demonstrated the possibility to infer low conducting regions, as sometimes witnessed for scarred regions, but future studies could apply the algorithm to analyze different pathologies, such as fibrosis.

With special care for registration, \TheMethod\ could also be combined with high-resolution 3D imaging, such as CT or MRI, to improve anatomical accuracy. An interesting question is whether \TheMethod\ could also be applied to ventricular activation. A major difference between ventricular and atrial activation is transmural propagation in the former, which is not accessible by contact mapping. Moreover, endocardial activation in the ventricles of healthy subjects, due to the Purkinje network, is extremely complex and may overshadow myocardial propagation. Under specific pathological conditions, such as ventricular tachycardia or bundle branch block, myocardial activation becomes relevant and heterogeneity in conduction of potential interest, justifying the applicability of \TheMethod. While it is true that no transmural data would be available, it is also known that 
fibers in ventricles follow a peculiar pattern in the transmural direction with low inter-patient variability~\cite{streeter69}. Such \emph{prior knowledge} may be used in the inverse procedure by appropriately changing the regularization term.
In a recent work \cite{grandits_inverse_2020}, we actually applied an inverse method similar to \TheMethod\ in the ventricles by using epicardial data, obtaining convincing results also in the transmural direction.

In light of the presented results, we believe that \TheMethod\ can assist future medical interventions by estimating
cardiac conduction properties more robustly and help in identifying ablation sites, as well
as in better understanding atrial and ventricular conduction pathways.
%
%
%
\bibliographystyle{res/llncs2e/splncs04}
\bibliography{MICCAI2020,inverse}




\end{document}